\documentclass[righttag,12pt]{amsart}
\usepackage{amssymb}
\usepackage{euscript}

 \let\cal\mathcal

\newcommand{\rest}[2]{\left.#1\right|_{#2}} 
\newtheorem{prop}{Proposition}[section]
\newtheorem{theorem}[prop]{Theorem}
\newtheorem{cor}[prop]{Corollary}
\newtheorem{lem}[prop]{Lemma}
 
\numberwithin{equation}{section}

\theoremstyle{definition}
\newtheorem{defn}[prop]{Definition}
\newtheorem{ack}{Acknowledgments}

\theoremstyle{remark}
\newtheorem{rem}[prop]{Remark}
\newcommand{\bbN}{{\mathbb{N}}}
 \newcommand{\bbR}{{\mathbb{R}}}
\newcommand{\bbT}{{\mathbb{T}}}

\newcommand{\la}{\lambda}

 \newcommand{\Om}{\Omega}

\newcommand{\Span}{\operatorname{span}}

\renewcommand{\Im}{\operatorname{Im}}
\newcommand{\sgn}{\operatorname{sgn}}
\renewcommand{\Re}{\operatorname{Re}}

\renewcommand{\le}{\ensuremath{\leqslant}}
\renewcommand{\ge}{\ensuremath{\geqslant}}
\renewcommand{\leq}{\ensuremath{\leqslant}}
\renewcommand{\geq}{\ensuremath{\geqslant}}

\newcommand{\lb}{\label}

\newcommand{\lra}{\longrightarrow}

\newcommand{\ONTO}{\buildrel {\mbox{\small onto}}\over \longrightarrow}

\begin{document}

\title{On contractive projections in Hardy spaces}

\author{Florence Lancien}
\address{D\'{e}partement de Math\'{e}matiques\\ Universit\'{e} de Franche-Comt\'{e}\\
16 Route de Gray, 25030 Besan\c{c}on, France}

\email{Florence.Lancien@math.univ-fcomte.fr}

\author{Beata Randrianantoanina$^*$}\thanks{$^*$Participant, NSF Workshop
in Linear Analysis and Probability, Texas A\&M University}
\address{Department of Mathematics and Statistics \\ Miami University
\\Oxford, OH 45056, USA}

\email{randrib@muohio.edu}

\author{Eric Ricard}
\address{D\'{e}partement de Math\'{e}matiques\\ Universit\'{e} de Franche-Comt\'{e}\\
16 Route de Gray, 25030 Besan\c{c}on, France}

\email{Eric.Ricard@math.univ-fcomte.fr}



\begin{abstract}
We prove a  conjecture of Wojtaszczyk
  that   for $1\leq p<\infty$, $p\neq 2$, $H_p(\bbT)$ does not
admit any norm one  projections   with dimension of the range finite
and bigger than 1. This implies in particular   that  for
$1\leq p<\infty$, $p\ne 2$, $H_p$ does not admit a Schauder basis with constant one.
\end{abstract}


\maketitle

\section{Introduction}

The study of norm one projections  and their ranges
(1-complemented subspaces) has been an important topic of the
isometric Banach space theory since the inception of the field.
Contractive projections were also investigated from the
approximation theory
 point of view, as part of  the study of minimal projections, i.e. projections onto
 the given subspace with the smallest possible norm (cf. \cite{ChP70,OL90}). They
 are also closely related to
the metric projections or nearest point mappings, and they are a
natural extension of the notion of orthogonal projections from
Hilbert spaces to general Banach spaces. Despite a great amount of
work on contractive projections in various function spaces (cf.
the survey \cite{survey}), very little is known about them in
spaces of analytic functions.

To the best of our knowledge, the only results about contractively
complemented subspaces of spaces of analytic functions are due to
Wojtaszczyk \cite{W79}, who proved that $H_\infty$ and the disc
algebra $\mathcal{A}$ have ``few'' contractively complemented
finite dimensional subspaces.

\begin{theorem} \cite[Theorem~4.3]{W79}
(a) If $P:{\mathcal{A}}\lra \mathcal{A}$ is a norm one  projection
with finite-dimensional range of dimension greater than one, then
$\Im P^*$ is contained in the set of measures $\mu$, which are
singular with respect to the Lebesgue measure $\la$.

(b) If $P:{H_\infty}\lra H_\infty$ is a norm one  projection with
finite-dimensional range of dimension greater than one, then $\Im
P^*$ is contained in the space of functionals on $H^\infty$ such
that every norm preserving extension to $L^\infty(\la)$ is
 singular with respect to the Lebesgue measure
$\la$.
\end{theorem}

In particular it follows that the space $L_1(\la)/H_1^0(\la)$
whose dual is $ H_\infty$ does not admit any norm one projection
with $1<\dim\Im P <\infty$ \cite[Corollary~4.4]{W79}, and that
${\mathcal{A}}$ does not have a monotone basis
\cite[Theorem~4.6]{W79}.

Wojtaszczyk asked about the situation in $H_p$, for $1\le p<\infty$,
$p\ne 2$. Namely, he asked (\cite{Whbk} for $p=1$, \cite{W-p} for $p>1$) whether there exist any norm
one  projections on $H_p$ with $1<\dim\Im P <\infty$. He
conjectured that the answer is no, which would imply in particular
that $H_p$ does not admit a 1-unconditional basis extending the
result of  \cite{W82} that $H_1$ does not admit a 1-unconditional
basis.

In this paper we prove that Wojtaszczyk's conjecture is correct,
i.e. that indeed, for $1\leq p<\infty$ and $p\neq 2$, $H_p$ does not
admit any norm one  projections   with $1<\dim\Im P <\infty$. In Section~4
we note that $H_p(\bbT^2)$ does have norm one projections with nontrivial
finite dimensional range and we
 include some remarks about the form of norm one projections
on $H_p$ with infinite dimensional range.

We note that it follows from an earlier work of Neuwirth
\cite{N98} that $H_p$ does not admit a 1-unconditional basis for
$1\leq p<\infty$, $p\ne 2$. Actually, Theorem 7.5 therein shows that
$H_p$ does not have the unconditional metric approximation property,
a  notion weaker than the existence of a 1-unconditional basis.

\begin{ack} The research on this paper started during a sabbatical
visit of the second named author at the D\'{e}partement de
Math\'{e}matiques, Universit\'{e}
 de Franche-Comt\'{e} in Besan\c{c}on, France. She wishes to thank the Department and
all members of
the Functional Analysis Group
 for their hospitality during that visit.

We thank the referee for valuable comments.
\end{ack}

\section{Preliminaries}

\begin{defn}
\lb{defJ} Let $X$ be a Banach space.  We define a {\it duality
map} $J$ from $X$ into subsets of $X^*$ by the condition that $f
\in J(x)\subset X^*$ if and only if $ \|f \|_{X^{*}} = \| X \|_X$
and $\left< f, x \right> = \|x \|^2_X$.
\end{defn}

We note that if $X$ is a  strictly convex Banach space, then for
all $x\in X$, the set $J(x)$ contains exactly one functional. In
this case we will consider $J$ as a map from $X$ to $X^*$. If, in
addition, $X$ is reflexive and $X^*$ is strictly convex, then
$J^*:X^*\to X$, and $J^*=J^{-1}$.

Calvert \cite{Cal75} proved an important characterization of
contractively complemented subspaces of reflexive Banach spaces in
terms of the duality map $J$.

\begin{theorem}\lb{calvert} \cite{Cal75}
Let $X$ be a strictly convex reflexive Banach space with strictly
convex dual.  Then a closed linear subspace $Y$ of $X$ is a
$1$-complemented in $X$ if and only if $J(Y)$ is a linear subspace
of $X^*$.
\end{theorem}

This theorem is one of the most important tools in the theory of
norm one projections.  It has been proven independently  also in
\cite{AF92} (using different methods).

 Actually, we will need only a simple version of this last result (cf.
\cite[Proposition~1.1]{AF92}):

\begin{lem}\label{dim}
 Let $X$ be a Banach space, such that for any $x\in X$, $J(x)$ consists
of exactly one point. Then, if $Y$ is finite dimensional and
$1$-complemented in $X$, then $J(Y)$ linearly spans a finite dimensional
subspace of $X^*$ of dimension less or equal than $\dim Y$.
\end{lem}

\begin{proof}
Let $P$ be the contractive projection onto $Y$. For any $y\in Y$:
$$\|y\|^2=\langle y,J(y)\rangle=\langle P(y),J(y)\rangle=\langle y,P^*(J(y))
\rangle.$$
But as $\|P^*(J(y))\|\le \|J(y)\|$, using the hypothesis, we get that
 $P^*(J(y))=J(y)$. So,
$\Span J(Y)$ is contained in the range of $P^*$.
\end{proof}

%
%

We will need basic facts about duality in $H_p$, that we will
recall now. More details can be found in \cite{Du} or \cite{Ho}.

By identifying each function in $H_p$ of the unit disk, with the
function defined by its boundary values on the torus, one can see
$H_p$ as a subspace on $L_p(\mathbb{T})$ in the following way:
$$H_p=H_p(\mathbb{T})=\{ f\in L_p(\mathbb{T}),  \  \hat{f}(-n)=\int_0^{2\pi} f(e^{i\theta})
e^{in\theta} d\theta=0, \forall n\in \mathbb{N}^*\}.$$

For $1\leq p<\infty$ the dual of $H_p$ is isometrically isomorphic
to the quotient  $L_q/ H_{q,0}$, where $1/p+1/q=1$ and
$H_{q,0}=\{f\in H_p, \ f(0)=0\}.$ More precisely, each functional
$\phi$ in $H_p^*$ can be defined by a function $g$ in $L_q$ as
follows: $\phi(f)={\frac 1  {2\pi}}\int_0^{2\pi} f(e^{i\theta})
g(e^{i\theta}) d\theta$, for $f\in H_p$, and two functions $g$ and
$h$ in $L_q$ represent the same functional if and only if $g-h$ is
in $H_{q,0}$. This very elementary fact will be used in the
sequel.

\begin{lem}
 For $1\leq p<\infty$ and each $f\in H_p$,
there exists a unique element $J(f)$ in $H_p^*$ such that $
\|J(f) \|_{H_p^{*}} = \| f \|_{H_p}$ and $\left< J(f), f \right> =
\|f \|^2_{H_p}$.
\end{lem}
\begin{proof}
 As for $1<p<\infty$, $L_p$ spaces are strictly convex and uniformly smooth,
$H_p$ as well as $H_p^*$ are strictly convex. This gives the result in the
reflexive case.

For $p=1$, since $f$ is analytic and has full support, we know that
there is only one functional $g\in L_\infty$ such that
$ \|g \|_{L_\infty} = \| f \|_{H_1}$ and $\left< g, f \right> =
\|f \|^2_{H_1}$. Let $h\in J(f)\subset L_\infty/H_{\infty,0}$. Then
by Theorem 2 in \cite{Pa70}, there is a lifting $\tilde h$ of $h$ in
$L_\infty$ such that $\|\tilde h\|=\|h\|=\|f\|$. But, then
$$ \langle \tilde h, f \rangle =\left< h, f \right>=\|f \|^2_{H_1}.$$
So, we must have $\tilde h=g$ and $h$ is unique.
\end{proof}

Consequently, we can apply Lemma \ref{dim} to $H_p$ spaces for
$1\leq p <\infty$.

We point out that the map $J:H_p\to H_p^*$  is
non linear but satisfies $J(af)=|a|J(f)$ if $a \in \mathbb{C}$ and
$f\in H_p$.
For $f\in H_p$ we will denote
$$f^N=|f|^{p-2}\overline {f}=|f|^{p-1} \sgn \overline {f},$$ where
$\sgn (re^{i\theta})=e^{i\theta}$. This function will be very
useful because it is, up to a multiplicative constant, a
representant of $J(f)$ in $H_p^*$.

\section{Main results}

\begin{theorem}\lb{main}
Let $1\leq p <\infty$, $p\neq 2$.  Then the Hardy space $H_p$ contains
no $1$-complemented subspaces of finite dimension larger than one.
\end{theorem}

\begin{cor}
Let $1\leq p <\infty$, $p\neq 2$.  Then the Hardy space $H_p$ does not
have a Schauder basis with constant 1.
\end{cor}

\begin{proof}[Proof of Theorem~\ref{main}]

Suppose for contradiction, that there exists a finite-dimensional
1-com\-ple\-men\-ted subspace $Y$ in $H_p$ with $\dim Y = d\ge 2$.
There will be two different arguments depending on whether $p$ is
an even integer or not.

\medskip\noindent
{(a) case $p$ is not an even integer :}

Let us denote $f_1, \cdots, f_d$ linearly independent functions in
$H_p$ such that $Y= \Span \{f_1, \cdots, f_d\}$. By preceding
remarks, for each $f\in Y$, $\displaystyle |f|^{p-2}
\overline{f}={\frac {|f|^{p}}  f}$ is a representant of an element of
$J(Y).$ Thus
$$J(Y)=\{\frac {|f|^{p}}  f  \oplus H_{q,0}, \ f\in Y \} =
\{\frac{|\sum_{i=1}^d \alpha_i f_i|^{p}} {\sum_{i=1}^d \alpha_i
f_i} \oplus H_{q,0}, \ \alpha_i \in \mathbb{C}, \ i=1,\cdots,d
\}.$$
 Since $Y$ is
1-com\-ple\-men\-ted in $H_p$ with $\dim Y = d\ge 2$,  by Lemma~\ref{dim},
$J(Y)$ is included in a
linear subspace of $L_q/ H_{q,0}$.
We define the following complex vector subspace of $L_1/H_{1,0}$
by
$$Z=\Span_\mathbb{C} \{{\frac {|f|^{p}}  f }f_i \oplus H_{1,0}, \
f\in Y, \ i=1, \cdots, d \}.$$ Since, by Lemma~\ref{dim}, $\dim \Span J(Y)\le d$,
we conclude that the space $Z$ has dimension at most $d^2$ as a
complex vector space and at most $2d^2$ as a real vector space.

We now consider the two functions $f_1$ and $f_2$ which are
linearly independent by assumption and we define $W$ and $W'$ as the
following real vector span
\begin{eqnarray*}
W'=\Span _\mathbb{R} \{|\alpha_1
f_1+\alpha_2 f_2|^{p}\oplus H_{1,0},\ \alpha_1, \alpha_2 \in
\mathbb{C}\}\subset L_1/H_{1,0}\\
W=\Span _\mathbb{R} \{|\alpha_1
f_1+\alpha_2 f_2|^{p},\ \alpha_1, \alpha_2 \in
\mathbb{C}\}\subset L_1.
\end{eqnarray*}
$W'$ is contained in $Z$, indeed if $f=\sum_{i=1}^2
\alpha_i f_i$, then $|f|^{p}=\sum_{j=1}^2 \alpha_j \frac {|f|^{p}}
f  f_j$ is a representant of an element in $Z$.
So $W'$ is finite dimensional. In turn, it implies that $W$ also has
finite dimension. Indeed, there is a natural map $i:W\to W'$, which sends
a function to its equivalence class, and it is injective. This follows from
the fact that we are dealing with real functions; a function $f$
in the kernel of
$i$ is in $H_{1,0}$, but then $f$ is  analytic and real, so $f$ is constant with
null average, it means that $f=0$.

We next choose $\Omega$ a subset of the torus of  non zero measure,
such that there exists $a>0$ with $a<|f_i|<1/a$ on $\Omega$ for
$i=1,2$ and we define
 $$\tilde{W}=\{\frac 1 {|f_1|^{p}} \rest h \Omega , \ h\in W\}.$$
As well as $Y$, $Z$ and $W$, this space is of finite dimension. We
are going to show that this leads to a contradiction.

For this purpose we consider the function $\Psi:\mathbb{C} \to
\tilde{W}$ defined by $\displaystyle\Psi(z)=|1+z \frac {f_2}{
f_1}\big|_\Omega|^p$. Write $p=2s$, since $f_2/f_1$ is bounded and bounded
away from zero on $\Omega$, we have the following power expansion
for $\Psi$ valid for $|z|$ small enough (below all expressions  $f_2/f_1$
are considered restricted to $\Om$):
$$\left|1+z \frac{f_2}{f_1}\right|^p=\left|\big(1+z \frac{f_2}
{f_1}\big)^s\right|^2=
\sum_{k,l\in\mathbb{N}} \bigl( _k^s \bigr) \bigl( _l^s \bigr)
\Big[\frac{f_2}{f_1}\Big]^k \overline{\Big[\frac{f_2}{f_1}\Big]^l}z^k\overline{z}^l.$$
Write $z=re^{i\theta}$ for $r\in \bbR$ small enough and integrate in
$\theta\in [0,2\pi]$. We obtain a function $\Phi:I\to \tilde W$, defined
in a neighbourghood of 0, such that
$$\Phi(r)=\sum_{k\geq 0} \bigl( _k^s \bigr)^2
\big|\frac{f_2}{f_1}\big|^{2k} r^{2k}$$

 As $\Phi$ takes its values in $\tilde W$, all its coefficients
in its power expansion are  also in $\tilde W$.
Since $p$ is not
an even integer, the binomial coefficient $\bigl( _k^s\bigr)$ is
never zero. One deduces that
$|\frac{f_2}{f_1}|^{2k}$ belongs to $\tilde{W}$ for all $k\in
\mathbb{N}$.

Since $\tilde{W}$ is a finite dimensional real vector space,
there must exists $K\in\mathbb{N}$ and  $\alpha_1,
\cdots,\alpha_K$ non all zero real coefficients such that
$$\sum_{k=1}^K \alpha_k \big|\frac{f_2}{f_1}\big|^{2k}= 0.$$
It  implies that $|\frac{f_2}{f_1}|$ takes values
in the finite set of the roots of the polynomial $\sum_{k=1}^K
\alpha_k X^{2k}$. The set $\Omega$ being of measure non zero, we
conclude that there exists a subset $\Omega_0$ of $\Omega$ and a
constant $C_0$ with $\mu(\Omega_0)>0$ and $|\frac{f_2}{f_1}|=C_0$
on $\Omega_0$.

If we take $(f_1, f_1+f_2)$ instead of $(f_1,f_2)$, we obtain
similarly the existence of a subset $\Omega_1$ of $\Omega_0$ and a
constant $C_1$ with $\mu(\Omega_1)>0$ and $|1+\frac{f_2}{f_1}|=C_1$
on $\Omega_1$. Finally on $\Omega_1$ we have, $|\frac{f_2}{f_1}|=C_0$
and $|1+\frac{f_2}{f_1}|=C_1$. Since on $\mathbb{C}$ the
equations $|z|=c_0$ and $|1+z|=c_1$ define two circles, this shows
that on $\Omega_1$, $f_1/f_2$ takes its values in the intersection
which consists of at most 2 points. By analyticity this implies
that $f_1=\lambda f_2$ on $\Omega_1$, since $\mu(\Omega_1)\neq 0$.
Therefore it follows that $f_1$ and $f_2$ are linearly dependent, which
contradicts the definition of $f_1$ and $f_2$.

\begin{rem} The referee has pointed out to us that, alternatively, one can use the
Rudin-Plotkin Equimeasurability Theorem \cite{Ru76,KKHandbook} to show that finite dimensionality
of $W$ implies that $f_1/f_2$ assumes only finite number of values. A sketch of an
argument is as follows:

If $\dim W<\infty$, then also $\dim V<\infty$, where  $V=\Span \{|1+(z(f_1/f_2)\big|_\Om)|^p \, :\, z\in C\}$.
Then we consider  space $V$ in $L_1$ with respect to the
measure $|f_1|^pdt$. Take any $h \in L_\infty=(L_1)^*$ which annihilates
$V$ and write $h=h_+ -h_-$ (decomposition into positive and negative
part).
We get
$\int |1+zf_1/f_2|^p(h_+ |f_1|^p)dt= \int |1+zf_1/f_2|^p(h_+|f_1|^p)dt$
for all complex $z$. From the Equimeasurability Theorem we get that
$f_1/f_2$ is equimeasurable with respect to both densities in brackets.
This implies that for any Borel bounded function $\Psi$ on $\mathbb{C}$ we have
$\int \Psi(f_1/f_2) h|f_1|^p dt =0$. Since this holds for any $h$ which
annihilates $V$ we infer that $\Psi(f_1/f_2)\in V$ for each $\Psi$. But
all those compositions can be finite dimensional only if $f_1/f_2$ assumes
only finite numbers of values.
\end{rem}

\medskip\noindent
{(b) case $p$ is an even integer :}

We denote $P$ a contractive projection on $H_p$ such that $Y$ is
the range of $P$. For fixed $h\in Y$ and $\phi \in H_p$, we have
for all $z\in \mathbb{C}$:
$$\|h+z\phi\|_p^p \geq \|h+zP(\phi)\|_p^p.$$
The function $\Psi$ defined by $\Psi(z)=\|h+z\phi\|_p^p -
\|h+zP(\phi)\|_p^p$ is a polynomial in $z$ and $\bar{z}$ which satisfies
$\Psi(0)=0$ and $\Psi(z)\geq 0$ for all $z$, thus $\Psi'(0)=0$.
This yields, after restricting $\Psi$ to real $z$ only, and differentiating
with respect to $z$,
$$\Re \left[\int_\mathbb{T}|h|^{p-2} {\overline h}(\phi-P(\phi))\right]=0.$$
Since this holds for any $h$ in $Y$, by multiplying $h$ by
complex constants of modulus one we obtain that
\begin{equation}\label{main1}\ \int_\mathbb{T}|h|^{p-2}\overline
{h}(\phi-P(\phi))=0,\ \forall h\in Y,\  \forall\phi\in H_p.
\end{equation}

We now pick $f$ and $g$ two functions in $Y$
and define $\Theta$ by
$$\Theta(z,z')=\int_\mathbb{T}|z'f+zg|^{p-2}\overline{(z'f+zg)}
(\phi-P(\phi)).$$
By \ref{main1} we have that $\Theta(z',z)=0$ for all
$z,z'\in\mathbb{C}$. Moreover, since $p$ is an even integer, $\Theta$
is a polynomial in $z$ and $z'$. By expanding this polynomial, the coefficient of
$\overline z |z'|^{p-2}$ gives that
$$\int_\mathbb{T}|f|^{p-2}\overline
{g}(\phi-P(\phi))=0,\  \forall\phi\in H_p.$$

Consequently, for any $f\in Y$
\begin{equation}\label{main2} \int_\mathbb{T}|f|^{p-2}\overline{h}
(\phi-P(\phi))=0,\ \forall h\in Y,\  \forall\phi\in H_p.
\end{equation}

We now introduce the scalar product $<.,.>_f$ defined  on $H_p$ by
$$<h,k>_f=\int_\mathbb{T}|f|^{p-2}\overline {h} k.$$
It is a scalar product since $f$ has full support
as it is  analytic. We
denote by $H_f$ the space $H_p$ equipped with the scalar product
$<.,.>_f$ and by $I_f$ the canonical injection from $H_p$ to
$H_f$. With these new notations, \ref{main2} can be written
$$<h,\phi-P(\phi)>_f=0, \ \forall h\in Y,\  \forall\phi\in H_p.$$
In particular for $h=P(\phi)$ we obtain
$$<P(\phi),\phi-P(\phi)>_f=0, \  \forall\phi\in H_p.$$
Thus $$\|P(\phi)\|_f^2=<\phi,P(\phi)>_f \leq
\|\phi\|_f\|P(\phi)\|_f, \  \forall\phi\in H_p.$$
This means that
$P_f$ defined by $P(\phi)=P_f(\phi)$ on $H_p$ is an orthogonal
projection on $H_f$ which coincides with $P$.

 From now on, we consider two
linearly independent functions $f$ and $g$ in $Y$.
We have  another scalar product $<.,.>_g$, a space
$H_g$ and $P_g$ an orthogonal projection on $H_g$. We next
consider $h_1,\cdots ,h_d$ a basis of $Y$ which is orthogonal
simultaneously for $<.,.>_f$ and $<.,.>_g$, this is possible because
$Y$ is assumed to be finite dimensional. Indeed, if $(e_i)_{i=1}^d$ is a $<.,.>_f-$orthonormal basis,
and  $A$ is the matrix of $<.,.>_g$
in that basis, then $A$ is   self adjoint and positive, and thus  there
exists a  $<.,.>_f-$orthonormal basis  which diagonalizes $A$, and thus is
$<.,.>_g-$orthogonal. Then, as $P_f$ and $P_g$ are
orthogonal projections
$$P_f=\sum_{i=1}^d \frac{<.,h_i>_f h_i}{\|h_i\|_f},$$
$$P_g=\sum_{i=1}^d \frac{<.,h_i>_g h_i}{\|h_i\|_g}.$$
Since for all $\phi\in H_p$ we have $P_f(\phi)=P_g(\phi)=P(\phi),$ we
obtain that $$\frac{<\phi,h_1>_f h_1} {\|h_1\|_f}=\frac{<\phi,h_1>_g
h_1}{\|h_1\|_g}, \  \forall\phi\in H_p.$$
If we denote
$\alpha=\frac 1 {\|h_1\|_f}$ and $\beta=\frac 1{\|h_1\|_g}$, this
yields
$$\int_\mathbb{T}(\alpha|f|^{p-2}-\beta|g|^{p-2})\overline {h_1}\phi=0,\
\forall\phi\in H_p.$$ This means that
$(\alpha|f|^{p-2}-\beta|g|^{p-2})\overline {h_1}=0$ in $H_p^*$
thus there exists $h\in H_{q,0}$ so that
$$(\alpha|f|^{p-2}-\beta|g|^{p-2})|h_1|^2=h h_1.$$ Again the fact
that the left hand-side is real and the analyticity of the right
hand-side imply that they are both zero.

Finally we obtain that $|f|^{p-2}=\lambda|g|^{p-2}$ for some
positive constant $\lambda$. This reasoning holds for any two
functions $f$ and $g$ linearly independent in $Y$. If we take
$f+g$ and $g$ instead of $f$ and $g$, we obtain that
$|f+g|^{p-2}=\nu|g|^ {p-2}$ for some positive constant $\nu$. But
we have already explained that this cannot happen unless $f$ and
$g$ are linearly dependent which contradicts our hypothesis.
\end{proof}

\section{Comments}

The above result is purely isometric. By that, we mean that
for any $d>1$ and $1\le p\leq \infty$ and $\epsilon>0$, there are
$d$-dimensional subspaces
of $H_p$ which are $1+\epsilon$ complemented. This is an easy consequence
of Szeg\"o theorem as it provides functions in $H_p$ which have almost disjoint
support.

\smallskip

One can also wonder if Theorem \ref{main1} has extensions to other domains
than the torus $\bbT$. As soon as there are more variables, the situation
becomes quite different :

\begin{prop}
For all $1\le p \leq \infty$ and any $d>1$, there are $1$-complemented $d$
dimensional subspaces in $H_p(\bbT^2)$
\end{prop}

\begin{proof}
 This is a just a transference argument. Let $d>1$ be fixed and for each
function $f(z_1,z_2)$ on $\bbT^2$, define
$$T(f)(z_1,z_2)=\int_{T} f(zz_1,z z_2) \overline z^{d+1} dz.$$
 It is not hard to see that $T$ ,for any $p$,  is  a contractive
projection from $L_p$ onto the subspace  spanned by functions $z_1^kz_2^l$ with
$k+l=d-1$.
\end{proof}

 The same kind on argument carries on for $H^p(S_n)$ where $S_n$ is the unit
sphere in $\mathbb{C}^{n+1}$, for any $n\geq1$. As an easy consequence, $H_p(\bbT)$
and $H_p(\bbT^2)$ are not isometric.

\bigskip

We will now describe some infinite dimensional 1-complemented
subspaces of $H_p$. For $\Lambda$ a subset of $\mathbb{N}$, we
define
$$H_{p,\Lambda}=\{ f\in H_p(\mathbb{T}),  \  \hat{f}(n)=0,\  \forall n\in \mathbb{N}
\smallsetminus \Lambda \}.$$ It is well known that if $\Lambda$ is
a coset in $\mathbb{N}$ then $H_{p,\Lambda}$ is 1-complemented in
$H_p$ , in that case the projection $P:H_p \ONTO H_{p,\Lambda}$ is
the Fourier multiplier defined by the characteristic function of
$\Lambda$.

More generally, for $\psi$ an inner function with $\psi(0)=0$, one
can consider $H_p^\psi$ the closure in $H^p$ of the algebra
generated by $\psi$. J.A. Ball \cite{Ba75} proved that there
exists an expectation operator $P_\psi$ which is a norm one
projection from $H_p$ onto $H_p^\psi$.

From this result we can deduce generalizations of the spaces
$H_\Lambda$ above. For $\psi$ an inner function with $\psi(0)=0$
and $\Lambda=\{qk+r, n\in\mathbb{N}\}$ where $q\in\mathbb{N}$ and
$r\in\{0,..., k-1\}$, we can  consider
$$H_{p,\Lambda}^\psi \textrm{ the closure in }H^p\textrm{ of the set }\{\psi^n, \
n\in\Lambda\}.$$ It is easy to see that $H_{p,\Lambda}^\psi$ is a
$1$-complemented subspace of $H_p$ via the projection defined for
$f\in H_p$ by $\displaystyle P_{\psi,\Lambda}(f)=\psi^{r-q}
P_{\psi^q}({f \psi^{q-r}})$ if $r\neq 0$ and $\displaystyle
P_{\psi,\Lambda}(f)=P_{\psi^q}(f)$ if $r=0$. (Note that
$H_{p,\Lambda}$ corresponds to the case $\psi(z)=z$.)

Moreover, it is clear that if $X$ is any Banach space, $T$ is an
isometry from $X$ onto $X$, and $P\ONTO Y$ is a norm one
projection, then $Q=TPT^{-1}: X\ONTO T(Y)$ is also a norm one
projection, and thus $T(Y)$ is 1-complemented in $X$. Isometries
of $H_p$ have been described by Forelli \cite{F64} who proved the
following.

\begin{theorem} \cite[Theorem 2]{F64}
Suppose that $p\ne 2$ and $T$ is a linear isometry of $H_p$ onto
$H_p$. Then
\begin{equation}\lb{F}
Tf=b\left(\frac{d\phi}{dz}\right)^{\frac1p}f(\phi),
\end{equation}
where $\phi$ is a conformal map of the unit disk onto itself and
$b$ is a unimodular complex number. Conversely, \eqref{F} defines
a linear isometry of of $H_p$ onto $H_p$.
\end{theorem}

If we denote $T_\phi$ an isometry of this form we obtain new norm
one projections by composition with the preceding projections
$P_{\psi,\Lambda}$:
$$Q=T_\phi P_{\psi,\Lambda} T_\phi^{-1}: H_p\ONTO T_\phi(H_{p,\Lambda}^\psi)$$
is also a contractive projection, and thus
$T_\phi(H_{p,\Lambda}^\psi)$ is 1-complemented in $H_p$.


We do not know whether every 1-complemented subspace of $H_p$ has
to be  isometric to $H_{p,\Lambda}^\psi$, for some coset $\Lambda$
in $\bbN$ and some inner function $\psi$ vanishing at the origin.



\def\polhk#1{\setbox0=\hbox{#1}{\ooalign{\hidewidth
  \lower1.5ex\hbox{`}\hidewidth\crcr\unhbox0}}} \def\cprime{$'$}

\end{document}